\documentclass[12pt]{article}

\usepackage{graphicx,amsmath,amssymb,amsthm,amsfonts,latexsym,tikz,hyperref,subfig}
\usetikzlibrary{decorations.pathmorphing}
\usepackage[hmargin=.8in,vmargin=1in]{geometry}

\DeclareMathOperator{\Lrm}{Lrm} \DeclareMathOperator{\lrm}{lrm}
 \DeclareMathOperator{\HT}{ht}
 \DeclareMathOperator{\spea}{spea}
\DeclareMathOperator{\npea}{npea} 
\DeclareMathOperator{\stun}{stun} \DeclareMathOperator{\fix}{fix}
\DeclareMathOperator{\exc}{exc}

\newcommand{\fp}{{\mathfrak p}}
\newcommand{\fI}{{\mathfrak I}}
\newcommand{\Up}{\Upsilon}

\newcommand{\M}{{\mathcal{M}}}

\DeclareMathOperator{\area}{area}

\DeclareMathOperator{\pos}{pos} \DeclareMathOperator{\val}{val}

\newtheorem{thm}{Theorem}[section]
\newtheorem{prop}[thm]{Proposition}
\newtheorem{cor}[thm]{Corollary}
\newtheorem{lem}[thm]{Lemma}
\newtheorem{conj}[thm]{Conjecture}
\newtheorem{exa}[thm]{Example}

\newtheorem{defn}[thm]{Definition}

\newcommand{\ben}{\begin{enumerate}}
\newcommand{\een}{\end{enumerate}}
\newcommand{\ble}{\begin{lem}}
\newcommand{\ele}{\end{lem}}
\newcommand{\bth}{\begin{thm}}
\renewcommand{\eth}{\end{thm}}
\newcommand{\bpr}{\begin{prop}}
\newcommand{\epr}{\end{prop}}
\newcommand{\bco}{\begin{cor}}
\newcommand{\eco}{\end{cor}}
\newcommand{\bcon}{\begin{conj}}
\newcommand{\econ}{\end{conj}}
\newcommand{\bde}{\begin{defn}}
\newcommand{\ede}{\end{defn}}
\newcommand{\bex}{\begin{exa}}
\newcommand{\eex}{\end{exa}}
\newcommand{\barr}{\begin{array}}
\newcommand{\earr}{\end{array}}
\newcommand{\btab}{\begin{tabular}}
\newcommand{\etab}{\end{tabular}}
\newcommand{\beq}{\begin{equation}}
\newcommand{\eeq}{\end{equation}}
\newcommand{\bea}{\begin{eqnarray*}}
\newcommand{\eea}{\end{eqnarray*}}
\newcommand{\bal}{\begin{align*}}
\newcommand{\bce}{\begin{center}}
\newcommand{\ece}{\end{center}}
\newcommand{\bpi}{\begin{picture}}
\newcommand{\epi}{\end{picture}}
\newcommand{\bpp}{\begin{picture}}
\newcommand{\epp}{\end{picture}}
\newcommand{\bfi}{\begin{figure} \begin{center}}
\newcommand{\efi}{\end{center} \end{figure}}
\newcommand{\bprf}{\begin{proof}}
\newcommand{\eprf}{\end{proof}\medskip}

\newcommand{\bsl}{\begin{slide}{}}
\newcommand{\esl}{\end{slide}}
\newcommand{\bfr}{\begin{frame}}
\newcommand{\efr}{\end{frame}}

\newcommand{\pf}{{\bf Proof.\hspace{5pt}}}
\newcommand{\prf}{{\noindent\bf Proof.\hspace{5pt}}}
\newcommand{\hqed}{\hfill \qed}
\newcommand{\hqedm}{\hfill \qed \medskip}

\newcommand{\hs}[1]{\hspace{#1}}
\newcommand{\hso}[1]{\hspace{-1pt}}

\newcommand{\qmq}[1]{\quad\mbox{#1}\quad}

\newcommand{\emp}{\emptyset}

\newcommand{\sbe}{\subseteq}

\newcommand{\case}[4]{\left\{\barr{ll}#1&\mbox{#2}\\#3&\mbox{#4}\earr\right.}

\newcommand{\flf}[2]{\left\lfloor\frac{#1}{#2}\right\rfloor}
\newcommand{\cef}[2]{\left\lceil\frac{#1}{#2}\right\rceil}

\def\<{\langle}
\def\>{\rangle}

\newcommand{\ree}[1]{(\ref{#1})}

\newcommand{\ra}{\rightarrow}

\newcommand{\al}{\alpha}
\newcommand{\be}{\beta}

\newcommand{\de}{\delta}
\newcommand{\ep}{\epsilon}

\newcommand{\la}{\lambda}

\newcommand{\si}{\sigma}

\newcommand{\Ga}{\Gamma}
\newcommand{\De}{\Delta}

\newcommand{\1}{{\bf 1}}

\newcommand{\bbN}{{\mathbb N}}

\newcommand{\bbZ}{{\mathbb Z}}
\newcommand{\cA}{{\cal A}}

\newcommand{\cD}{{\cal D}}

\newcommand{\cH}{{\cal H}}

\newcommand{\cM}{{\cal M}}

\newcommand{\cP}{{\cal P}}

\newcommand{\cR}{{\cal R}}

\newcommand{\fS}{{\mathfrak S}}

\DeclareMathOperator{\Av}{Av}

\DeclareMathOperator{\col}{col} 
 \DeclareMathOperator{\des}{des}
\DeclareMathOperator{\notdes}{\overline{des}} 
\DeclareMathOperator{\Des}{Des}

 \DeclareMathOperator{\inv}{inv}

\DeclareMathOperator{\maj}{maj}

 \DeclareMathOperator{\st}{st}

 \DeclareMathOperator{\wt}{wt}

\begin{document}
\title{Inversion polynomials for $321$-avoiding permutations
}
\author{Szu-En Cheng\\[-5pt]
\small Department of Applied Mathematics, National University of Kaohsiung\\[-5pt]
\small Kaohsiung 811, Taiwan, ROC, {\tt chengszu@nuk.edu.tw}\\[5pt]
Sergi Elizalde\\[-5pt]
\small Department of Mathematics, Dartmouth College\\[-5pt]
\small Hanover, NH 03755-3551, USA, {\tt sergi.elizalde@dartmouth.edu}\\[5pt]
Anisse Kasraoui\\[-5pt]
\small Fakult\"at f\"ur Mathematik, Universit\"at Wien\\[-5pt]
\small Nordbergstra{\ss}e 15, A-1090 Vienna, Austria, {\tt anisse.kasraoui@univie.ac.at}\\[5pt]
Bruce E. Sagan\\[-5pt]
\small Department of Mathematics, Michigan State University,\\[-5pt]
\small East Lansing, MI 48824-1027, USA, {\tt sagan@math.msu.edu} }

\date{\today\\[10pt]
    \begin{flushleft}
    \small Key Words:  Catalan number, continued fraction, Dyck path, generating function, pattern avoidance, permutation, inversion number,  major index, Motzkin path, polyomino, $q$-analogue.
                                           \\[5pt]
    \small AMS subject classification (2010):
    Primary 05A05;
    Secondary 05A10, 05A15, 05A19, 11A55.
    \end{flushleft}}

\maketitle

\begin{abstract}
We prove a generalization of a conjecture of Dokos, Dwyer, Johnson,
Sagan, and Selsor giving a recursion for the inversion polynomial of
$321$-avoiding permutations.  We also answer a question they posed
about finding a recursive formulas for the major index polynomial of
$321$-avoiding permutations.  Other properties of these polynomials
are investigated as well.  Our tools include Dyck and 2-Motzkin
paths, polyominoes, and continued fractions.
\end{abstract}

\section{Introduction}
\label{i}

The main motivation for this paper is a conjecture of Dokos, Dwyer,
Johnson, Sagan, and Selsor~\cite{ddjss:pps} about inversion
polynomials for $321$-avoiding permutations which we will prove in
generalized form.  We also answer a question they posed by giving a
recursive formula for the analogous major index polynomials.   We
first introduce some basic definitions and notation about pattern
avoidance and permutation statistics.

Call two sequences of distinct integers $\pi=a_1\ldots a_k$ and
$\si=b_1\ldots b_k$ {\em order isomorphic} whenever $a_i<a_j$ if and
only if $b_i<b_j$ for all $i,j$.  Let $\fS_n$ denote the symmetric
group of permutations of $[n]\stackrel{\rm def}{=}\{1,\ldots,n\}$.
Say that {\em $\si\in\fS_n$ contains $\pi\in\fS_k$ as a pattern} if
there is a subsequence $\si'$ of $\si$ order isomorphic to $\pi$.
If $\si$ contains no such subsequence then we say $\si$ {\em avoids}
$\pi$ and write $\Av_n(\pi)$ for the set of such $\si\in\fS_n$.

Let $\bbZ$ and $\bbN$ denote the integers and nonnegative integers,
respectively.  A {\em statistic} on $\fS_n$ is a function
$\st:\fS_n\ra\bbN$.  One then has the corresponding generating
function
$$
f_n^{\st}=\sum_{\si\in\fS_n} q^{\st\si}.
$$
Two of the most ubiquitous statistics for $\si=b_1\ldots b_n$ are
the {\em inversion number}
$$
\inv\si=\#\{(i,j)\ |\ \text{$i<j$ and $b_i>b_j$}\}
$$
where the hash sign denotes cardinality, and the {\em major index}
$$
\maj\si=\sum_{b_i>b_{i+1}} i.
$$

In~\cite{ss:mp}, Sagan and Savage proposed combining the study of
pattern avoidance and permutations statistics by considering
generating functions of the form 
\beq \label{F^st}
F^{\st}(\pi)=\sum_{\si\in\Av_n(\pi)} q^{\st\si} 
\eeq 
for any pattern
$\pi$ and statistic $\st$.  Dokos et al.~\cite{ddjss:pps} were the
first to carry out an extensive study of these generating functions
for the $\inv$ and $\maj$ statistics. We note that when $\st=\inv$ and $\pi=132$ we recover a $q$-analogue of the Catalan numbers studied by Carlitz and Riordan~\cite{cr:tel}.
Work on the statistics
counting fixed points and excedances has been done by
Elizalde~\cite{eli:fpe, eli:mpa}, Elizalde and
Deutsch~\cite{ed:sub}, and Elizalde and Pak~\cite{ep:brr}.

Our primary motivation was to prove a conjecture of Dokos et al.\
concerning the inversion polynomial for $321$-avoiding permutations.
In fact, we will prove a stronger version which also keeps track of
left-right maxima. Call $a_i$ in $\pi=a_1\ldots a_n$ a {\em
left-right maximum (value)} if $a_i=\max\{a_1,\ldots,a_i\}$.  We let
$$
\Lrm \pi =\{a_i\ |\ \text{$a_i$ is a left-right maximum}\}
$$
and $\lrm\pi=\#\Lrm\pi$.   Consider the generating function \beq
\label{I321:def} I_n(q,t)=\sum_{\si\in\Av_n(321)} q^{\inv\si}
t^{\lrm\si}. \eeq Note that since $\#\Av_n(321)=C_n$, the $n$th
Catalan number, this polynomial is a $q,t$-analogue of $C_n$.  Our
main result is a recursion for $I_n(q,t)$.  The case $t=1$ was a
conjecture of Dokos et al. \bth \label{I321:rr} For $n\ge1$,
$$
I_n(q,t)=tI_{n-1}(q,t)+\sum_{k=0}^{n-2} q^{k+1} I_k(q,t)
I_{n-k-1}(q,t).
$$
\eth The rest of this paper is structured as follows.    In the next
section we will give a direct bijective proof of
Theorem~\ref{I321:rr} using $2$-Motzkin paths.  The following two
sections will explore related ideas involving Dyck paths, including
a combinatorial proof of a formula of  F\"urlinger and
Hofbauer~\cite{fh:qcn} and two new statistics which are closely
related to $\inv$.  Sections~\ref{ptp} and~\ref{mip} are devoted to
polyominoes.  First, we give a second proof of Theorem~\ref{I321:rr}
using work of Cheng, Eu, and Fu~\cite{cef:acp}.  Next we derive 
recursions for a major index analogue, $M_n(q,t)$, of~\ree{I321:def},
thus answering a question posed by Dokos et al.\ in their paper.  In
Section~\ref{sum}, we explore further properties of $M_n(q,t)$,
including symmetry, unimodality, and its modulo $2$ behavior.    The
final two sections are concerned with continued fractions.  We begin
by giving a third proof of Theorem~\ref{I321:rr} using  continued
fractions and comparing it with a result of
Krattenthaler~\cite{kra:prp}. 
In fact, this demonstration is even more general since it also keeps track of the number of fixed points.
 And in  Section~\ref{rse} we reprove
and then generalize a theorem of Simion and Schmidt~\cite{ss:rp}
concerning the signed-enumeration of $321$-avoiding permutations.

\section{A proof of Theorem~\ref{I321:rr} using 2-Motzkin paths}
\label{ptm}

Our first proof of Theorem~\ref{I321:rr} will use 2-Motzkin paths.
Let $U$ (up), $D$ (down), and $L$ (level) denote vectors in $\bbZ^2$
with coordinates $(1,1)$, $(1,-1)$, and $(1,0)$, respectively.  A
{\em Motzkin path of length $n$}, $M=s_1\ldots s_n$,  is a lattice
path where each step $s_i$ is $U$, $D$, or $L$ and which begins at
the origin, ends on the $x$-axis, and never goes below $y=0$.  A
{\em $2$-Motzkin path} is a Motzkin path where each level step has
been colored in one of two colors which we will denote by $L_0$ and
$L_1$.
We will let $\cM_n^{(2)}$ denote the set of $2$-Motzkin paths of
length $n$.

It will be useful to have two vectors to keep track of the values
and positions of left-right maxima.  If $\si=b_1\ldots b_n\in\fS_n$
then let
$$
\val\si =(v_1,\dots,v_n)
$$
where
$$
v_i=\case{1}{if $i$ is a left-right maximum of $\si$,}{0}{else.}
$$
Also define
$$
\pos\si=(p_1,\dots,p_n)
$$
where
$$
p_i=\case{1}{if $b_i$ is a left-right maximum of $\si$,}{0}{else.}
$$
By way of example, if $\si=361782495$ then we have
$\val\si=(0,0,1,0,0,1,1,1,1)$ and $\pos\si=(1,1,0,1,1,0,0,1,0)$.
Note that for any permutation $v_n=p_1=1$.

We will  need the following lemma which collects together some
results from the folklore of pattern avoidance.  Since they are easy
to prove, the demonstration will be omitted. \ble
\label{321:structure} Suppose $\si\in\fS_n$. \ben
\item[(a)]  We have $\si\in\Av_n(321)$ if and only if the elements of $[n]-\Lrm \si$ form an increasing subsequence of $\si$.
\item[(b)]  Suppose we are given $0$-$1$ vectors $v=(v_1,\ldots,v_n)$ and $p=(p_1,\ldots,p_n)$ with the same positive number of ones.
Then $v=\val\si$ and $p=\pos\si$ for some $\si\in\fS_n$ if and only
if, for every index $i$, $1\le i\le n$, the number of ones in
$p_1,\ldots,p_i$ is greater than the number in $v_1,\ldots,v_{i-1}$.
In this case, because of part (a), there is a unique such
$\si\in\Av_n(321)$. \hqedm \een \ele Note that in (b) the cases
where $i=1$ and $i=n$ force $p_1=1$ and $v_n=1$, respectively.

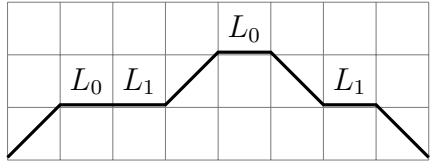
\begin{figure}
\begin{center}
\begin{tikzpicture}[scale=.7]
\draw[step=1,color=gray] (0,0) grid (8,3); \draw [very thick]
(0,0.05)--(1,1.05) --
(2,1.05)--(3,1.05)--(4,2.05)--(5,2.05)--(6,1.05)--(7,1.05)--(8,0.05);
\draw (1.5,1.5) node{$L_0$}; \draw (2.5,1.5) node{$L_1$}; \draw
(4.5,2.5) node{$L_0$}; \draw (6.5,1.5) node{$L_1$};
\end{tikzpicture}
\end{center}
\caption{The Motzkin path associated with $\si=361782495$}
\label{fig:motzkin}
\end{figure}

\noindent{\bf First proof of Theorem~\ref{I321:rr}.}\hs{5pt}  We
will construct a bijection $\mu:\Av_n(321)\ra\cM_{n-1}^{(2)}$ as
follows.  Given $\si\in\Av_n(321)$ with $\val\si=(v_1,\ldots,v_n)$
and $\pos\si=(p_1,\ldots,p_n$), we let $\mu(\si)=M=s_1\ldots
s_{n-1}$ where
$$
s_i=
\begin{cases}
U&\text{if $v_i=0$ and $p_{i+1}=1$,}\\
D&\text{if $v_i=1$ and $p_{i+1}=0$},\\
L_0&\text{if $v_i=p_{i+1}=0$,}\\
L_1&\text{if $v_i=p_{i+1}=1$.}
\end{cases}
$$
Continuing the example from the beginning of the section,
$\si=361782495$ would be mapped to the path in
Figure~\ref{fig:motzkin}.

We must first show that $\mu$ is well defined in that if
$M=\mu(\si)$ then $M$ ends on the $x$-axis and stays weakly above it
the rest of the time.  In other words, we want the number of $U$'s
in any prefix of $M$ to be at least as great as the number of $D$'s,
with equality at the finish.  This now follows from the definition
of the $s_i$ and the first two sentences of
Lemma~\ref{321:structure} (b).

We must also check that $\mu$ is a bijection.  The fact that it is
injective is an immediate consequence of the definition of the $s_i$
and the third sentence of Lemma~\ref{321:structure} (b).  Since
$\#\Av_n(321)=C_n=\#\cM_{n-1}^{(2)}$, we also have bijectivity.

If  $\mu(\si)=M$ then we claim that
\begin{eqnarray}
\lrm\si&=&\#U(M)+\#L_1(M)+1,\\
\label{invsi} \inv\si&=&\#D(M)+\#L_0(M)+\area M,
\end{eqnarray}
where $\area M$ is the area between $M$ and the $x$-axis, $U(M)$ is
the set of steps equal to $U$ in $M$, and similarly for the other
types of steps.  The first equation follows from the definition of
$M$ and the fact that $\lrm\si$ is the number of ones in $\pos\si$.
The  $+1$ is because $M$ has length $n-1$ and we always have $p_1=1$.

To prove the equation for $\inv$,  we will induct on $n$.  Note that
every $M\in \M_{n-1}^{(2)}$, where $n\ge2$, can be uniquely
decomposed in one of the following ways: \ben
\item[(i)] $M=L_0 N$, where $N\in \M_{n-2}^{(2)}$,
\item[(ii)] $M=L_1 N$, where $N\in \M_{n-2}^{(2)}$.
\item[(iii)] $M=U N D O$, where $N\in \M_{k-1}^{(2)}$ and $O\in \M_{n-k-2}^{(2)}$ for some $1\le k\le n-2$.
\een

Suppose~\ree{invsi} holds for $N$ in case (i), and suppose
$\mu(\pi)=N$.  If we have $\val\pi=(v_1,\ldots,v_{n-1})$ and
$\pos\pi=(p_1,\ldots,p_{n-1})$ then adding  $L_0$ forces the vectors
to change to $\val\si=(0,v_1,\ldots,v_{n-1})$ and
$\pos\si=(1,0,p_2,\ldots,p_{n-1})$.  It follows that if
$\pi=a_1\ldots a_{n-1}$ then $\si=(a_1+1)1(a_2+1)\ldots(a_{n-1}+1)$.
So both the left and right sides of~\ree{invsi} go up by one when
passing from $\pi$ to $\si$ and equality is preserved.  Similar
arguments shows that both sides stay the same in case~(ii), and both
go up by $k+1$ in case~(iii).  So equality is maintained in all
cases.

From what we have shown, it suffices to show that
$$
I'_n(q,t)=\sum_{P\in \M_{n-1}^{(2)}} q^{\#D(P)+\#L_0(P)+\area P}
t^{\#U(P)+\#L_1(P)+1}
$$
satisfies the recurrence in Theorem~\ref{I321:rr}.  Considering the
three cases above, in~(i) we get a contribution of $q I'_{n-1}(q,t)$
to $I'_n(q,t)$ which corresponds to the $k=0$ term of the sum.  Similarly, case~(ii) contributes  $t
I'_{n-1}(q,t)$.  Finally, in~(iii)
 the piece $UND$ contributes $tq^{k+1} I_{k}(q,t)$ since when lifting $N$ the area is increased by $k$, and both $\#U(N)$ and $\#D(N)$ are increased by one. Also $O$ contributes $I_{n-k-1}(q,t)/t$ since there is a $+1$ in the exponent of $t$ for both $I_{k}(q,t)$ and $I_{n-k-1}(q,t)$, but we only want one such.
Combining these contributions proves the recursion.\hqedm

\section{Dyck paths and an equation of F\"urlinger and Hofbauer}
\label{dpe}

In this section we will prove a formula of F\"urlinger and
Hofbauer~\cite[equation (5.5)]{fh:qcn} which is closely related to
Theorem~\ref{I321:rr}.  In fact, we will show in Section~\ref{ptp}
that this equation can be used to prove our main theorem.  Our proof
of the F\"urlinger-Hofbauer result will be combinatorial using Dyck
paths, whereas the one given in~\cite{fh:qcn} is by algebraic
manipulation of generating functions.  Our proof has the interesting
feature that it uses a nonstandard decomposition of Dyck paths which will also be useful in the next section.

Let $P=s_1\ldots s_{2n}$ be a Dyck path of semilength $n$ and let
$\cD_n$ be the set of all such $P$.  We will freely go back and
forth between three standard interpretations of such paths.  In the
first, $P$ consists of $n$ $U$-steps and $n$ $D$-steps starting at
the origin and staying weakly above the $x$-axis.  It the second,
there are $n$ north steps, $N=(0,1)$, and $n$ east steps, $E=(1,0)$,
beginning at the origin and staying weakly above the line $y=x$.  In
the last, we have $n$ zeros and $n$ ones with the number of zeros in
any prefix of $P$ being at least as great as the number of ones.  In
this last interpretation, we can apply all the usual permutation
statistics defined in the same way as they were when there were no
repetitions.  In particular, we will need the {\em descent set} of
$P$
$$
\Des P =\{i\ |\ a_i > a_{i+1}\}
$$
and the {\em descent number} $\des P = \#\Des P$.  A descent of $P$
as a bit string corresponds to a {\em valley} of $P$ in the first
interpretation, i.e., a factor of the form $DU$.  We will also need
the dual notion of a {\em peak}, which is a factor $UD$.

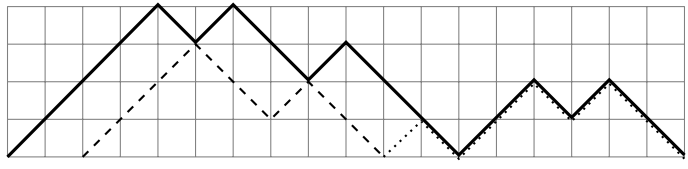
\begin{figure}
\begin{center}
\begin{tikzpicture}[scale=.5]
\tikzstyle{every node}=[font=\tiny] \draw[step=1,color=gray] (0,0)
grid (18,4); \draw [very thick] (0,0)--(4,4.05) -- ++(1, -1) --
++(1, 1) -- ++(2, -2) -- ++(1, 1) -- ++(3, -3) -- ++(2, 2) -- ++(1,
-1) -- ++(1, 1) -- ++(2, -2); \draw [thick,dashed] (2,0)--(5,3) --
++(2, -2) -- ++(1, 1) -- ++(2, -2); \draw [thick,dotted]
(10,0)--(11, 0.95) -- ++(1, -1) -- ++(2, 2) -- ++(1, -1) -- ++(1, 1)
-- ++(2, -2);
\end{tikzpicture}
\end{center}
\caption{The decomposition of Dyck paths} \label{fig:decomposition}
\end{figure}

We also need to define one of the analogues of the Catalan numbers
 studied by F\"urlinger and Hofbauer.  Given a Dyck path
$P=s_1\ldots s_{2n}$ we let $|P|_0$ and $|P|_1$ be the number of
zeros and number of ones in  $P$, respectively.  More generally we
will write $|w|_A$ for the number of occurrences of $A$ in the word
$w$ for any $A$ and $w$. Now let \beq \label{fp} \fp_i(w)=s_1\ldots
s_i \eeq
 be $w$'s prefix of length $i$.  Define
\bea
\al(P)&=&\sum_{i\in\Des P} |\fp_i(P)|_0,\\
\be(P)&=&\sum_{i\in\Des P} |\fp_i(P)|_1. \eea Note that
$\al(P)+\be(P)=\maj P$.  Now consider the generating function \beq
\label{C_n(t)} C_n(t)=C_n(a,b;t)=\sum_{P\in\cD_n} a^{\al(P)}
b^{\be(P)} t^{\des P}. \eeq
\bth[F\"urlinger and
Hofbauer~\cite{fh:qcn}] \label{fh:thm} We have \beq \label{C_n(t)}
C_n(t)=C_{n-1}(abt) + bt\sum_{k=0}^{n-2} a^{k+1} C_k(abt)
C_{n-k-1}((ab)^{k+1}t). \eeq \eth \prf We will first define a
bijection $\de:\biguplus_{k=0}^{n-1}\cD_k\times\cD_{n-k-1}\ra\cD_n$.
Given two Dyck paths
$$
Q=U^{a_1}D^{b_1}U^{a_2}D^{b_2}\dots U^{a_s}D^{b_s}\in\cD_k \qmq{and}
R=U^{c_1}D^{d_1}U^{c_2}D^{d_2}\dots U^{c_t}D^{d_t}\in\cD_{n-k-1}
$$
where all  exponents are positive, we will combine them
to  create a Dyck path $P=\de(Q,R)\in\cD_n$  as follows.  (When $Q$
is empty, the same definition works with the convention that
$a_1=b_1=0$).   There are two cases: \ben
\item If $R=\emp$, then
$$
P=U^{a_1+1} D^{b_1+1} U^{a_2}D^{b_2}\dots U^{a_s}D^{b_s}.
$$
\item If $R\neq\emp$, then
$$
P=U^{a_1+1}DU^{a_2}D^{b_1}U^{a_3}D^{b_2}\dots U^{a_s}D^{b_{s-1}}
U^{c_1}D^{b_s+d_1}U^{c_2}D^{d_2}\dots U^{c_t}D^{d_t}.
$$
\een For example,
$$
\de(U^3D^2UD^2,UDU^2DUD^2)=U^4DUD^2UD^3U^2DUD^2,
$$
as illustrated in Figure~\ref{fig:decomposition}.  The path $P$ is
given by the solid lines while $Q$ is dashed  and $R$ is dotted, with $Q$
being  shifted to start at $(2,0)$ and  $R$ concatenated directly
after, starting at $(10,0)$.  Note that this puts the peaks of $Q$ at
exactly the same position as the valleys of the first part of $P$,
and makes $R$ and $P$ coincide after the first peak of $R$

To show that $\de$ is bijective, we construct its inverse. Suppose
that
$$
P=U^{i_1}D^{j_1}U^{i_2}D^{j_2}\dots U^{i_k}D^{j_k}.
$$
Again, there are two cases for computing $\de^{-1}(P)=(Q,R)$. \ben
\item If $j_1\ge2$, then
$$
Q=U^{i_1-1}D^{j_1-1}U^{i_2}D^{j_2}\dots U^{i_k}D^{j_k} \qmq{and}
R=\emp.
$$
\item If $j_1=1$, let $s$ be the smallest index such that
$i_1+i_2+\dots+i_s\le j_2+j_3+\dots+j_{s+1}$ (note that $s<k$), and
let $\ep=j_2+j_3+\dots+j_{s+1}-(i_1+i_2+\dots+i_s)$. Then
$$
Q=U^{i_1-1}D^{j_2}U^{i_2}D^{j_3}\dots U^{i_s}D^{j_{s+1}-\ep-1}
\qmq{and} R=U^{i_{s+1}}D^{\ep+1}U^{i_{s+2}}D^{j_{s+2}}\dots
U^{i_k}D^{i_k}.
$$
\een The first case is easy to understand as you just shorten the
first peak of $P$ which had been lengthened by $\de$.
 For the second case, the reader may find it useful to  consult Figure~\ref{fig:decomposition} again.  As in defining $\de$, $Q$ is the Dyck path that starts at $(2,0)$ and has the peaks at the valleys of $P$.  At some point this will no longer be possible because $Q$ would have to go under the $x$-axis. $Q$ ends at the point on the $x$-axis just before it would be forced to go negative, and $R$ starts at that point. The path $R$ begins with up-steps until it hits $P$, and then coincides with $P$ for the rest of the path.

We now use our bijection to prove the theorem.  If $P\in\cD_n$ and
$\de^{-1}(P)=(Q,R)$, then the term $C_{n-1}(abt)$ corresponds to the
case when $R=\emp$, since then adding the extra $UD$ to the first
peak moves all the valleys over by one unit. When $R\neq\emp$,
suppose $Q\in\cD_k$ and $R\in\cD_{n-k-1}$. The factor
$C_{n-k-1}((ab)^{k+1}t)$ comes from the fact that all the valleys of
$R$ become valleys of $P$, each one having $k+1$ additional steps
$U$ and $D$ to their left. The term $C_k(abt)$ accounts for the fact
that when the valleys of $Q$ become valleys of $P$, they have  an
additional $U$ and $D$ inserted to their left, namely the first $D$
in $P$ and the $U$ preceding it, and the term $abt$ accounts for the
extra valley started by this $D$. Finally, the factor $a^k$ comes
from the fact that when the $k$ up-steps of $Q$ are put in $P$, they
are shifted one valley to the left of the corresponding valley in
$P$.  Thus  each of these steps  moves  the corresponding
valley one position to the right. \hqedm

The decomposition in Theorem~\ref{fh:thm} is based on the structure of the first peak of $P$.
If we consider the last peak of $P$, we will get the following equation which does not appear in \cite{fh:qcn}.

\bth We have
$$
C_n(t)=C_{n-1}(t) + a^{n-1}t\sum_{k=1}^{n-1} b^{k} C_k(t)
C_{n-k-1}((ab)^{k}t).
$$
\eth
\prf We define a second
bijection $\De:\biguplus_{k=0}^{n-1}\cD_k\times\cD_{n-k-1}\ra\cD_n$ as follows.
Given two Dyck paths
$$
R=U^{c_1}D^{d_1}U^{c_2}D^{d_2}\dots U^{c_t}D^{d_t}\in\cD_{k} \qmq{and}
Q=U^{a_1}D^{b_1}U^{a_2}D^{b_2}\dots U^{a_s}D^{b_s}\in\cD_{n-k-1}
$$
\ben
\item If $R=\emp$, then
$$
P=U^{a_1} D^{b_1} U^{a_2}D^{b_2}\dots U^{a_s +1}D^{b_s +1}.
$$
\item If $R\neq\emp$, then
$$
P=U^{c_1}D^{d_1}U^{c_2}D^{d_2}\dots U^{c_t+a_1}D^{d_t}U^{a_2}D^{b_1}U^{a_3}D^{b_2}\dots U^{a_s}D^{b_{s-1}}UD^{b_s +1}.
$$
\een
Now using similar arguments to those in the proof of Theorem~\ref{fh:thm}, we  get the desired result. \hqedm

\section{The sumpeaks and sumtunnels statistics}
\label{sss}

In this section we discuss another pair of (new) statistics on
Dyck paths, which are closely related to the $\inv$ statistic on
$321$-avoiding permutations.  In fact, these two statistics are
equidistributed over $\fS_n$, as we will show.  First we will need
some definitions and notation.

Let $\npea P$ denote the number of peaks of $P$, and note that $\npea P=\des P+1$
for any non-empty Dyck path $P$.
Also define the {\em height} of a peak $p=UD=s_is_{i+1}$ in a Dyck path
$P=s_1\ldots s_{2n}$ to be
\begin{eqnarray}
\HT(p)
\label{ht}
=|\fp_i(P)|_U - |\fp_i(P)|_D,
\end{eqnarray}
where $\fp_i$ is as in equation~\ree{fp}. In the interpretation of
Dyck paths using steps $U$ and $D$, $\HT(p)$ is the $y$-coordinate
of the highest point of $p$. In Figure~\ref{fig:decomposition}, the
peaks of $P$ have heights $4, 4, 3, 2,$ and $2$.

\bfi \label{Gamma:fig}
\begin{center}
\begin{tikzpicture}[scale=.5]
\tikzstyle{every node}=[font=\tiny] \draw[step=1,color=gray] (0,0)
grid (9,9); \draw [very thick]
(0,0)--(0,3)--(1,3)--(1,4)--(3,4)--(3,6)--(6,6)--(6,9)--(9,9); \draw
(0,0)--(9,9); \fill (0.5,2.5) circle (5pt); \fill (1.5,3.5) circle
(5pt); \fill (2.5,0.5) circle (5pt); \fill (3.5,5.5) circle (5pt);
\fill (4.5,1.5) circle (5pt); \fill (5.5,4.5) circle (5pt); \fill
(6.5,8.5) circle (5pt); \fill (7.5,6.5) circle (5pt); \fill
(8.5,7.5) circle (5pt);
\end{tikzpicture}
\end{center}\
\caption{The bijection $\Ga$} \efi

To make a connection with permutations $\pi=a_1\ldots a_n$, we
consider the {\em diagram} of $\pi$, which is a square grid with a
dot in column $j$ at height $a_j$ for all $j\in[n]$.  Figure 2
displays the permutation $\pi=341625978\in\Av_9(321)$.  We will now
describe a bijection $\Ga:\Av_n(321)\ra\cD_n$ which appeared
in~\cite{eli:fpe} (where it is denoted by $\psi_3$), and in a
slightly different form in~\cite{kra:prp}.  Imagine a light shining
from the northwest of the diagram of $\pi$ so that each dot casts a
shadow with sides parallel to the axes.  Consider the lattice path
$P$ formed by the boundary of the union of these shadows.  (This is
the same procedure as used by Viennot~\cite{vie:fgc} in his
geometric version of the Robinson-Schensted correspondence.)  Define
$\Ga(\pi)=P$.  Again, Figure 2 illustrates the process.  Using
Lemma~\ref{321:structure}, one can prove that for any permutation
$\pi\in\fS_n$, the path $P$ will stay above $y=x$ and so be a Dyck
path.  Lemma~\ref{321:structure} also shows that if one restricts to
$\pi\in\Av_n(321)$, then this map becomes a bijection.  The $\inv$
and $\lrm$ statistics on $\pi$ translate nicely under $\Ga$.
\bpr\label{prop:Ga} If  $\Ga(\pi)=P$, then \ben
\item[(a)] $\lrm \pi=\npea P$,
\item[(b)] $\inv \pi=\sum_{p} (\HT(p)-1)$ where the sum is over all peaks $p$ of $P$.
\een \epr \prf For (a), just note that every left-right maximum of
$\pi$ is associated with a peak $p=NE$ of $P$ consisting of two of
the edges of the square containing the maximum.  And this
correspondence is clearly reversible.

For (b), since  $\pi=a_1\ldots a_n$ avoids $321$, each inversion
$(i,j)$ of $\pi$ has the property that $a_i$ is a left-right maximum
but $a_j$ is not.  There is a two-to-one correspondence between
elements of $\pi$ and steps of $P$, where $a_i$ corresponds to the
pair $(N_i,E_i)$ which are the projections horizontally and
vertically onto $P$, respectively.  Now $(i,j)$ is an inversion if
and only if $N_j$ comes  before $N_i$ and $E_j$ comes after $E_i$.
Thus it follows from equation~\ree{ht} (with $U$ and $D$ replaced by
$N$ and $E$, respectively) that the element $a_i$ corresponding to a
peak $p$ causes inversions with exactly $\HT(p)-1$ elements $a_j$.
The $-1$ comes from the fact that $(i,i)$ is not an inversion.
Summing over all peaks completes the proof. \hqedm

We will show at the end of Section~\ref{mip} that the bijection
$\Ga$ can be used to give alternative proofs of
Theorems~\ref{I321:rr} and~\ref{thm:Mn}.

Because of its appearance in the previous proposition, we define a
new statistic {\em sumpeaks} on Dyck paths $P$ by
$$
\spea P =\sum_p (\HT(p)-1)
$$
where the sum is over all peaks $p$ of $P$.  There is another
statistic that we will now define which is equidistributed with
sumpeaks.

\begin{figure}
\begin{center}
\begin{tikzpicture}[scale=.5]
\tikzstyle{every node}=[font=\tiny] \draw[step=1,color=gray] (0,0)
grid (18,4); \draw [very thick] (0,0)--(4,4.05) -- ++(1, -1) --
++(1, 1) -- ++(2, -2) -- ++(1, 1) -- ++(3, -3) -- ++(2, 2) -- ++(1,
-1) -- ++(1, 1) -- ++(2, -2); \draw [thick,dashed]
(2,2.05)--(8,2.05)  (3,3.05)--(5,3.05)  (0,0.05)--(12,0.05)
(13,1.05)--(15,1.05);
\end{tikzpicture}
\end{center}
\caption{The tunnels of a Dyck path} \label{fig:tunnel}
\end{figure}
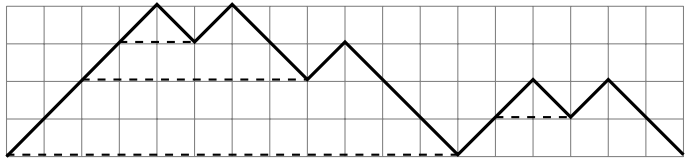

If $P$ is a Dyck path with steps $U,D$ and  $v=DU=s_js_{j+1}$ is a valley of
$P$ then its {\em height}, $\HT(v)$, is the $y$-coordinate of its
lowest point.  For each
valley, $v$, there  is a corresponding {\em tunnel}, which is the
factor $T=s_i\ldots s_j$ of $P$ where $s_i$ is the step after the
first intersection of $P$ with the line $y=\HT(v)$ to the left of
$s_j$.  The tunnels for the Dyck path in Figure~\ref{fig:tunnel} are
indicated with dashed lines.  In every tunnel, $j-i$ is an even
number, so we define the {\em sumtunnels} statistic to be
$$
\stun P = \sum_{T=s_i\ldots s_j} (j-i)/2
$$
where the sum is over all tunnels $T$ of $P$.  In
Figure~\ref{fig:tunnel}, we have $\stun P= (12+6+2+2)/2=11$.  It
turns out that the sumpeaks and sumtunnels statistics are
equidistributed over $\cD_n$.
\bth\label{thm:speastun} For any $n\ge1$,
$$
\sum_{P\in\cD_n} q^{\spea P} t^{\npea P} = \sum_{P\in\cD_n} q^{\stun
P} t^{n-\des P}.
$$
\eth 
\prf For $P\in\cD_n$, let $\notdes P=n-\des P$ for convenience.
Recall that $\des P$ is the number of valleys of $P$. It suffices to
define a bijection $h:\cD_n\ra\cD_n$ such that for any $P\in\cD_n$
we have \beq \label{spea} \spea P=\stun h(P) \quad\mbox{and}\quad
\npea P=\notdes h(P). \eeq Let $\de^{-1}(P)=(Q,R)$ where $\de$ is
the bijection of the Section~\ref{dpe}.  We inductively define $h$
by $h(\emp)=\emp$ and for $n\ge1$
$$
h(P)= \left\{ \barr{ll}
UDh(Q)&\text{if $R=\emp$,}\\
Uh(R)D&\text{if $Q=\emp$,}\\
Uh(Q)Dh(R)&\text{else.} \earr \right.
$$

To see that $h$ has an inverse, it suffices to check that given
$P'\in\cD_n$ we can tell which of the three cases above $P'$ must
fall into for $|P'|\ge4$. (When $|P'|=2$ then bijectivity is clear
since there is only one Dyck path of this length.) The first case
contains all $P'$ starting with a single $U$.  The second case
covers all $P'$ that are strictly above the $x$-axis between the
first and last lattice points which forces them to start with at
least two $U$'s.  And the last case contains those paths which start
with at least two $U$'s and intersect the $x$-axis before the final
vertex.

We now verify equation~\ree{spea} by induction on $n$.  It is easy
to verify for $n=1$.  For greater $n$, let
$$
P=U^{i_1}D^{j_1}U^{i_2}D^{j_2}\dots U^{i_k}D^{j_k},
$$
where $i_1,j_1,\dots,i_k,j_k$ are positive. We have three cases.

If $j_1\ge2$, then $Q=U^{i_1-1}D^{j_1-1}U^{i_2}D^{j_2}\dots
U^{i_k}D^{j_k}$ and $R=\emp$. So, comparing $P$ and $Q$ and  using
the induction hypothesis,
\begin{align*}
& \spea P =\spea Q + 1  =\stun h(Q) + 1 = \stun UDh(Q) = \stun h(P),\\
& \npea P=\npea Q=\notdes h(Q)=\notdes UDh(Q)=\notdes h(P).
\end{align*}

If $j_1=i_1=1$, then $Q=\emp$ and $R=U^{i_2}D^{j_2}\dots
U^{i_k}D^{j_k}$. Using similar reasoning to the first case,
\begin{align*}
& \spea P =\spea R =\stun h(R) = \stun Uh(R)D = \stun h(P),\\
& \npea P=\npea R+1=\notdes h(R)+1=\notdes Uh(R)D=\notdes h(P).
\end{align*}

If $j_1=1$ and $i_1 \ge2 $ then, keeping the notation in the
definition of $\de$,
$$
Q=U^{i_1-1}D^{j_2}U^{i_2}D^{j_3}\dots U^{i_s}D^{j_{s+1}-\ep-1}
\qmq{and} R=U^{i_{s+1}}D^{\ep+1}U^{i_{s+2}}D^{j_{s+2}}\dots
U^{i_k}D^{i_k}.
$$
The last peaks of $P$ coincide with all but the first peak of $R$.
The first peak of $R$ and the peaks of $Q$ are in bijection with the
rest of the peaks of $P$ where a peak of $P$ corresponds to the peak
of $Q$ (or the first peak of $R$) which is closest on its right.
Let $p_1,\ldots,p_{s+1}$ be these peaks of $P$, corresponding to peaks
$q_1,\ldots,q_s$ in $Q$ and $r_1\stackrel{\rm def}{=}q_{s+1}$ in
$R$.  Then
$$
\HT(p_k)= \left\{ \barr{ll}
\HT(q_1)+1&\text{if $k=1$,}\\
\HT(q_k)+j_k&\text{if $1<k\le s$,}\\
\HT(q_{s+1})+j_{s+1}-\ep-1&\text{if $k=s+1$.} \earr \right.
$$
Thus we have
$$
\sum_{k=1}^{s+1} \HT(p_k)=\sum_{k=1}^{s+1}
\HT(q_k)+j_2+\cdots+j_{s+1}-\ep=\sum_{k=1}^{s+1} \HT(q_k)+i_1+\cdots+i_s.
$$
Hence \bea
\spea P&=&\spea Q + \spea R + \sum_{k=1}^s i_k \\
&= &\stun h(Q) + \stun h(R) + \sum_{k=1}^s i_k\\
&=& \stun Uh(Q)Dh(R)\\
&=&\stun h(P), \eea where the third equality comes from the fact
that $Uh(Q)Dh(R)$ has exactly one more tunnel than the union of the
tunnels of $Q$ and $R$, namely the tunnel from the new $U$ to the
new $D$, and that tunnel has semilength $\sum_{k=1}^s i_k$.
Additionally,
$$
\npea P=\npea Q+\npea R=\notdes h(Q)+\notdes h(R)=\notdes
Uh(Q)Dh(R)=\notdes h(P),
$$
completing the proof.
\hqedm

\section{A proof of Theorem~\ref{I321:rr} using polyominoes}
\label{ptp}

In this section we will give a second proof of our main theorem
using Theorem~\ref{fh:thm}, another result of F\"urlinger and
Hofbauer, and polyominoes.  In particular, we will need a bijection
$\Up$ first defined by Cheng, Eu, and Fu~\cite{cef:acp} between
shortened polyominoes and $321$-avoiding permutations.  We first
need to define some terms.

\begin{figure}
\begin{center}
\begin{tikzpicture}[scale=0.8]
\tikzstyle{every node}=[font=\tiny] \draw [very thick]
(0,0)--(0,3)--(2,3)--(2,4)--(3,4)--(3,6)--(4,6); \draw [very thick]
(0,0)--(2,0)--(2,1)--(3,1)--(3,3)--(4,3)--(4,6); \draw (1,0)--(1,3);
\draw (0,1)--(2,1); \draw (2,1)--(2,3); \draw (0,2)--(3,2); \draw
(2,3)--(3,3); \draw (3,3)--(3,4); \draw (3,4)--(4,4); \draw
(3,5)--(4,5); \fill (0,0) circle (3pt) ; \fill (0,1) circle (2pt);
\fill (0,2) circle (2pt); \fill (0,3) circle (2pt); \fill (1,0)
circle (2pt); \fill (1,3) circle (2pt); \fill (2,0) circle (2pt);
\fill (2,1) circle (2pt); \fill (2,3) circle (2pt); \fill (2,4)
circle (2pt); \fill (3,1) circle (2pt); \fill (3,2) circle (2pt);
\fill (3,3) circle (2pt); \fill (3,4) circle (2pt); \fill (3,5)
circle (2pt); \fill (3,6) circle (2pt); \fill (4,3) circle (2pt);
\fill (4,4) circle (2pt); \fill (4,5) circle (2pt); \fill (4,6)
circle (3pt); \node at (2,-1) {\normalsize (a) A  polyomino in
$\cP_{10}$};
\end{tikzpicture}
\hspace{10pt}
\begin{tikzpicture}[scale=0.8]
\draw [very thick] (0,0)--(0,2)--(2,2)--(2,3)--(3,3)--(3,5)--(4,5);
\draw [very thick] (0,0)--(2,0)--(2,1)--(3,1)--(3,3)--(4,3)--(4,5);
\draw (1,0)--(1,2); \draw (0,1)--(2,1); \draw (2,1)--(2,2); \draw
(2,2)--(3,2); \draw (3,4)--(4,4); \fill (0,0) circle (2pt); \draw
(-.3,.5) node  {1}; \fill (0,1) circle (2pt); \draw (-.3,1.5) node
{2}; \fill (0,2) circle (2pt); \draw (.5,2.3) node  {3}; \fill (1,0)
circle (2pt); \draw (1.4,2.3) node  {4}; \fill (1,2) circle (2pt);
\draw (1.8,2.5) node  {5}; \fill (2,0) circle (2pt); \draw (2.4,3.3)
node  {6}; \fill (2,1) circle (2pt); \draw (2.8,3.5) node  {7};
\fill (2,2) circle (2pt); \draw (2.8,4.5) node  {8}; \fill (2,3)
circle (2pt); \draw (3.5,5.3) node  {9}; \fill (3,1) circle (2pt);
\draw (.5,-.3) node  {3}; \fill (3,2) circle (2pt); \draw (1.5,-.3)
node  {4}; \fill (3,3) circle (2pt); \draw (2.2,.4) node  {1}; \fill
(3,3) circle (2pt); \draw (2.6,.7) node  {6}; \fill (3,4) circle
(2pt); \draw (3.3,1.5) node  {2}; \fill (3,5) circle (2pt); \draw
(3.3,2.3) node  {5}; \fill (4,3) circle (2pt); \draw (3.6,2.7) node
{9}; \fill (4,4) circle (2pt); \draw (4.3,3.5) node  {7}; \fill
(4,5) circle (2pt); \draw (4.3,4.5) node  {8}; \node at (2,-1)
{\normalsize (b) A  polyomino in $\cH_9$};
\end{tikzpicture}
\end{center}
\caption{Parallelogram and shortened polyominoes} \label{poly:fig}
\end{figure}
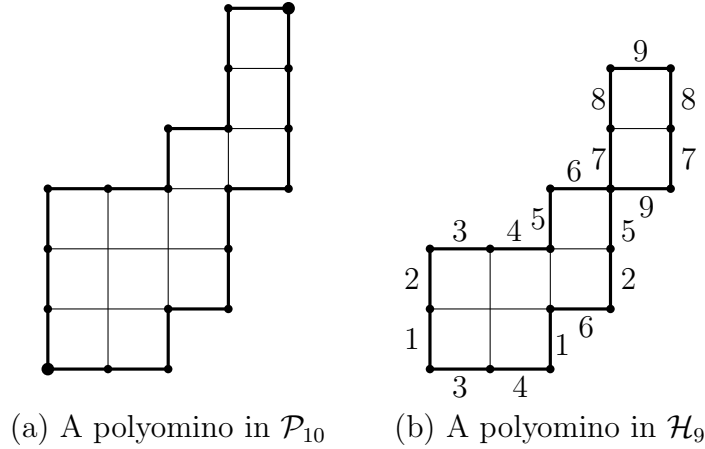

A {\em parallelogram polyomino} is a pair $(U,V)$ of lattice paths
using steps $N$ and $E$ such that
\begin{itemize}
\item  $U$ and $V$ begin at the same vertex and end at the same vertex, and
\item $U$ stays strictly above $V$ except at the beginning and end vertices.
\end{itemize}
In Figure~\ref{poly:fig}(a) we have $U=NNNEENENNE$ and
$V=EENENNENNN$.  Let $\cP_n$ denote the set of all parallelogram
polyominoes with $|U|=|V|=n$.  Note that if $U=s_1\ldots s_n$ and
$V=t_1\ldots t_n$ then $s_1=t_n=N$, and $s_n=t_1=E$. Define two
statistics \bea
\area(U,V)&=&\text{the area contained inside $(U,V)$},\\
\col(U,V)&=&\text{the number of columns spanned by $(U,V)$}. \eea
Returning to our example, $\area(U,V)=12$ and $\col(U,V)=4$.
Consider the generating function
$$
P_n(q,t)=\sum_{(U,V)\in\cP_n} q^{\area(U,V)}t^{\col(U,V)}.
$$
Another result of F\"urlinger and Hofbauer, which we state here
without proof, shows that this polynomial is closely related to
$C_n(a,b;t)$ as defined in equation~\ree{C_n(t)}. \bth[F\"urlinger
and Hofbauer~\cite{fh:qcn}] \label{fh:thm2} We have
$$
P_{n+1}(q,t)=q^nt C_n(q,q^{-1};t)
$$
for all $n\ge0$.\hqedm \eth

We will also need another type of polyomino.  Define a {\em
shortened polyomino} to be a pair $(P,Q)$ of $N,E$ lattice paths
satisfying
\begin{itemize}
\item  $P$ and $Q$ begin at the same vertex and end at the same vertex, and
\item $P$ stays weakly above $Q$ and the two paths can share $E$-steps but not $N$-steps.
\end{itemize}
Figure~\ref{poly:fig}(b) shows such a polyomino.  We denote the set
of shortened polyominoes with $|P|=|Q|=n$ by $\cH_n$.

We can now define the map $\Up:\cH_n\ra\Av_n(321)$.  Given
$(P,Q)\in\cH_n$, label the steps of $P$ with the numbers
$1,\ldots,n$ from south-west to north-east.  Each step of $P$ is
paired with the projection of that step onto $Q$.  Give each step of
$Q$  the same label as its pair.  Then reading the labels on $Q$
from south-east to north-west gives a permutation $\si=\Up(P,Q)$.
In Figure~\ref{poly:fig}(b), $\si=341625978$.  The next result
compares our statistics on $\cH_n$ and $\Av_n(321)$
\bth[Cheng-Eu-Fu~\cite{cef:acp}] \label{cef:thm} The map
$\Up:\cH_n\ra\Av_n(321)$ is a well-defined bijection such that if
$\Up(P,Q)=\si$ then \ben
\item[(a)] $\area(P,Q)=\inv\si$, and
\item[(b)] $\col(P,Q)=\lrm \si$.
\een \eth \prf The fact that $\Up$ is a well-defined bijection and
part (a) were proved in~\cite{cef:acp}, so we will only sketch the
main ideas here. If $\Up(P,Q)=\si$, then the left-right maxima of
$\si$ will label the $E$ steps of $Q$.  The positions of these
maxima in $\si$ are the same as their positions on $Q$.  Thus, as we
saw in Lemma~\ref{321:structure}(b), this data will determine a
unique $321$-avoiding permutation provided that the prefix condition
is satisfied.  And that condition is ensured by the second item in
the definition of a shortened polyomino.  Thus we have a bijection.

Now suppose $\si=a_1\ldots a_n$ and that we have an inversion
$a_i>a_j$ where $i<j$.  In that case $a_i$ and $a_j$ will label an
$E$-step and an $N$-step of $Q$, respectively, with the $N$-step
coming later on the path.  One can then  show that there will be a
square inside $(P,Q)$ due north of $a_i$ and due west of $a_j$
corresponding to the inversion.  This process is reversible, so
there is a bijection between inversions of $\si$ and squares inside
$(P,Q)$, proving part (a) of the theorem.  And part (b) follows from
the already-noticed fact that the left-right maxima of $\si$ are in
bijection with the $E$-steps of $Q$. \hqedm

The final ingredient is a simple bijection between $\cP_{n+1}$ and
$\cH_n$:  If $(U,V)\in\cP_{n+1}$ then contracting the first step of
$U$ and the last step of $V$ (both of which are $N$-steps) gives
$(P,Q)\in\cH_n$.  The polyomino in Figure~\ref{poly:fig}(b) is
gotten by shortening the one in~\ref{poly:fig}(a) in this manner.
If shortening $(U,V)$ gives $(P,Q)$ then we clearly have
\begin{eqnarray}
\label{area}
\area(U,V)&=&\area(P,Q)+\col(P,Q),\\
\col(U,V)&=&\col(P,Q). \label{col}
\end{eqnarray}

\noindent{\bf Second proof  of Theorem~\ref{I321:rr}.}
Theorem~\ref{cef:thm} together with equations~\ree{area}
and~\ree{col} give  $I_n(q,t)=P_{n+1}(q,t/q)$.  Combining this with
Theorem~\ref{fh:thm2} yields
$$
I_n(q,t)=\case{q^{n-1}t C_n(q,1/q;t/q)}{if $n\ge1$,}{1}{if $n=0$.}
$$
Now in Theorem~\ref{fh:thm} we replace $a,b,$ and $t$ by $q, 1/q$, and
$t/q$ , respectively.  Multiplying both sides by $q^{n-1}t$ and
rewriting everything in terms of the corresponding inversion
polynomials finishes the proof. \hqedm

\section{A major index polynomial recursion}
\label{mip}

In the paper of Dokos et al., they asked for a recursion for the
$321$-avoiding major index polynomial which is defined by
equation~\ree{F^st} with $\st=\maj$ and $\pi=321$.  The purpose of
this section is to give such a recurrence relation using
polyominoes.

Consider the polynomial \beq \label{M_n}
M_n(q,t)=\sum_{\si\in\Av_n(321)} q^{\maj\si}t^{\des\si}. \eeq Using
the description of the bijection $\Up:\cH_n\ra\Av_n(321)$ given in
the proof of Theorem~\ref{cef:thm}, it is clear that each descent of
$\si$ corresponds to a factor $EN$ of $Q$ where
$\Up^{-1}(\si)=(P,Q)$ and vice-versa.  Since the position of the
descent in $\si$ is the same as the position of the factor on $Q$,
we have $\Des \si=\Des Q$  where, as usual, $Q$ is identified with
the bit string obtained by replacing $N$ and $E$ by $0$ and $1$,
respectively.  It follows that $\des\si=\des Q$ and $\maj\si=\maj
Q$.  So we can rewrite~\ree{M_n} as
$$
M_n(q,t)=\sum_{(P,Q)\in\cH_n} q^{\maj Q}t^{\des Q}.
$$
We will need a lemma about what happens if we restrict this sum to
the parallelogram polyominoes $\cP_n\sbe\cH_n$. \ble We have
$$
\sum_{(P,Q)\in\cP_n} q^{\maj Q}t^{\des
Q}=M_{n-1}(q,t)+(q^{n-1}t-1)M_{n-2}(q,t).
$$
\ele \prf Let $(P',Q')\in\cH_{n-1}$ be obtained from $(P,Q)\in\cP_n$
by shortening.  Also let us  write $P=p_1\ldots p_{n-1}E$ and
$Q=q_1\ldots q_{n-1}N$.  We have two cases.

If $q_{n-1}=E$, then $n-1\in\Des Q$ which implies
$$
q^{\maj Q} t^{\des Q}=q^{n-1+\maj Q'}t^{1+\des Q'}=q^{n-1}tq^{\maj
Q'} t^{\des Q'}
$$
Furthermore, $P'$ and $Q'$ both end with an $E$ step and removal of
that common step leaves a polyomino in $\cH_{n-2}$.  It follows that
the $(P,Q)$ in this case contribute $q^{n-1}t M_{n-2}(q,t)$ to the
sum in the lemma.

For the second case  we have $q_{n-1}=N$.  It follows that
$(P',Q')\in\cH_{n-1}$ where the only restriction is that $Q'$ end
with a north step.  In other words, we want the generating function
for all polyominoes in $\cH_{n-1}$ except for those whose lower path
ends with an $E$ step (which must coincide with the last step of the
upper path which is always $E$).  This is clearly
$M_{n-1}(q,t)-M_{n-2}(q,t)$, and adding the contributions of the two
cases we are done. \hqedm

\bth\label{thm:Mn} For $n\ge1$ we have
$$
M_n(q,t)=M_{n-1}(q,qt)+\sum_{k=2}^{n}
\left[M_{k-1}(q,t)+(q^{k-1}t-1)M_{k-2}(q,t)\right]M_{n-k}(q,q^kt),
$$
and
$$
M_n(q,t)=M_{n-1}(q,t)+\sum_{k=0}^{n-2}  M_k(q,t)\left[M_{n-k-1}(q,q^kt)+(q^{n-1}t-1)M_{n-k-2}(q,q^kt)\right].
$$
\eth
\prf  To obtain the first equation, suppose $(P,Q)\in\cH_n$.  If both $P$ and $Q$ start with
an $E$ step then the generating function for such pairs is
$M_{n-1}(q,qt)$ since each descent of $Q$ is moved over one
position.

Since $Q$ always starts with an $E$-step, the only other possibility
is for $P$ to start with an $N$-step.  Let $z$ be the first point of
intersection of $P$ and $Q$ after their initial vertex.  Let $Q_0$
and $Q_1$ denote the portions of $Q$ before and after $z$,
respectively, and similarly for $P_0$ and $P_1$.  Let
$k=|Q_0|=|P_0|$.  But then $P_0$ and $Q_0$ do not intersect between
their initial point and $z$.  Thus $(P_0,Q_0)\in\cP_k$ and, from the
previous lemma,  the generating function for such pairs is the first
factor in the summation.

We also have  $|Q_1|=|P_1|=n-k$ and $(P_1,Q_1)\in\cH_{n-k}$.  Since
$Q_1$ is preceded by a path with $k$ steps, each of its descents
will be increased by $k$.   So the generating function for such
pairs is $M_{n-k}(q,q^kt)$.   Putting all the pieces together
results in the first formula in the statement of the theorem.

To obtain the second, merely replace $z$ in the proof  just given by the last point of intersection of $P$ and $Q$ before their final vertex.
 \hqedm

We end this section with the observation that there is a close connection between the maps $\Ga$ and $\Up$ defined in Sections~\ref{sss} and~\ref{ptp}, respectively.
Specifically, the inverse of
$\Ga:\Av_n(321)\ra\cD_n$ coincides
with the composition of the bijection from $\cD_n$ to $\cH_n$ used
in~\cite{fh:qcn} to prove Theorem~\ref{fh:thm2} with the bijection
$\Up:\cH_n\ra\Av_n(321)$.
So one can use $\Ga$ in place of $\Up$ in some of the proofs.
For example, the equation
$I_n(q,t)=q^{n-1}t C_n(q,1/q;t/q)$ which appears in the second proof
of Theorem~\ref{I321:rr} can also be obtained from $\Ga$ as
follows. If $\pi\in\Av_n(321)$ and $\Ga(\pi)=P$, then by
Proposition~\ref{prop:Ga} and the definitions in Section~\ref{dpe},
we have
\begin{align*}
& \inv(\pi)=\spea P=n+\alpha(P)-\beta(P)-\npea P=n-1+\alpha(P)-\beta(P)-\des P,\\
& \lrm(\pi)=\npea P=1+\des P.
\end{align*}
Thus,
$$I_n(q,t)=\sum_{P\in\cD_n} q^{\spea P} t^{\npea P}=q^{n-1}t C_n(q,1/q;t/q).$$

Theorem~\ref{thm:Mn} can also be proved using the bijection $\Ga$.
Note that each descent in a permutation $\pi\in\Av_n(321)$
corresponds to an occurrence of the string $NEE$ in the Dyck path
$\Ga(\pi)$. Thus the statistics $\des$ and $\maj$ in $\pi$
correspond to the number of occurrences and the sum of the $x$-coordinates after the first $E$ in each occurrence
 of $NEE$ in $\Ga(\pi)$, respectively. Using the standard
decomposition of a non-empty Dyck path as $P=NQER$ where $Q$ and
$R$ are Dyck paths, as well as its reversal, we can keep track of these two statistics to
obtain the recursions in Theorem~\ref{thm:Mn}.

\section{Symmetry, unimodality, and mod $2$ behavior of  $M_n(q,t)$}
\label{sum}

The coefficients of the polynomials $M_n(q,t)$ have various nice
properties which we now investigate. If $f(x)=\sum_k a_kx^k$ is a
polynomial in $x$ then we will use the notation
\bea
[x^k] f(x)&=&\mbox{coefficient of $x^k$ in $f(x)$}\\
&=&a_k.
\eea
Our main object of study in this section will be the polynomial
$$
A_{n,k}(q)=[t^k] M_n(q,t).
$$
In other words, $A_{n,k}(q)$ is the generating function for the
$\maj$ statistic over $\si\in\Av_n(321)$ having exactly $k$ descents. 

The first property which will concern us is symmetry. Consider a
polynomial
$$
f(x)=\sum_{i=r}^s a_i x^i
$$
where $a_r, a_s\neq0$.  Call $f(x)$ \emph{symmetric} if $a_i=a_j$
whenever $i+j=r+s$. \bth \label{sym} The polynomial $A_{n,k}(q)$ is
symmetric for all $n,k$. \eth \pf\ If $\si$ is counted by
$A_{n,k}(q)$ then $\des\si=k$.  Since $\si$ avoids $321$, it can not
have two consecutive descents and so the minimum value of $k$ is
$$
1+3+\cdots+(2k-1)=k^2
$$
and the maximum value is
$$
(n-1)+(n-3)+\cdots+(n-2k+1)=nk-k^2.
$$
So it suffices to show that for $0\le i\le nk$ we have
$a_i=a_{nk-i}$ where
$$
A_{n,k}(q)=\sum_i a_i q^i.
$$

Let $\cA_i$ be the permutations counted by $a_i$ and let $R_{\180}$
denote rotation of the diagram of $\si$ by $180$ degrees.  We claim
$R_{180}$ is a bijection between $\cA_i$ and $\cA_{nk-i}$ which will
complete the proof.  First of all, $R_{180}(321)=321$ and so
$\si$ avoids $321$ if and only if $R_{180}(\si)$ does so as well.  If $\si\in\cA_i$
then let $\Des\si=\{d_1,\ldots,d_k\}$ where $\sum_j d_j=i$.  It is
easy to see that $\Des R_{180}(\si)=\{n-d_1,\ldots,n-d_k\}$.  It
follows that $\maj R_{180}(\si)=nk-i$ and so $R_{180}(\si)\in
\cA_{nk-i}$.  Thus $R_{180}$ restricts to a well defined map from
$\cA_i$ to $\cA_{nk-i}$.  Since it is its own inverse, it is also a
bijection. \hqedm

Two other properties often studied for polynomials are unimodality
and log concavity.  The polynomial $f(x)=\sum_{i=0}^s a_i x^i$ is
{\em unimodal} if there is an index $r$ such that $a_0\le\ldots\le
a_r\ge\ldots\ge a_s$.  It is {\em log concave} if $a_i^2\ge a_{i-1}
a_{i+1}$ for all $0<i<s$.  If all the $a_i$ are positive, then log
concavity implies unimodality. \bcon The polynomial $A_{n,k}(q)$ is
unimodal for all $n,k$. \econ This conjecture has been checked by
computer for all $k< n\le 10$.  The corresponding conjecture for log
concavity is false, in particular, $A_{6,2}$ is not log concave.

The number-theoretic properties of the Catalan numbers have attracted
some interest.  Alter and Kubota~\cite{ak:ppp} determined the
highest power of a prime $p$ dividing $C_n$ using arithmetic means.
Deutsch and Sagan~\cite{ds:ccm}  gave a proof of this result using
group actions for the special case $p=2$.  Just considering parity,
one gets the nice result that $C_n$ is odd if and only if $n=2^m-1$
for some nonnegative integer $m$.  Dokos et al. proved the following
refinement of the ``if" direction of this statement. \bth[Dokos et
al.~\cite{ddjss:pps}] Suppose $n=2^m-1$ for some $m\ge0$.  Then
$$
\hspace*{130pt}[q^k]I_n(q,1)=\case{1}{if $k=0$}{\text{an even
integer}}{if $k\ge 1$.\hspace*{130pt}\qed}
$$
\eth In the same paper, the following statement was made as a
conjecture which has now been proved by Killpatrick.

\bth[Killpatrick~\cite{kil:pcc}] Suppose $n=2^m-1$ for some $m\ge0$.
Then
$$
\hspace*{130pt}[q^k]M_n(q,1)=\case{1}{if $k=0$}{\text{an even
integer}}{if $k\ge 1$.\hspace*{130pt}\qed}
$$
\eth

\bfi
\begin{tikzpicture}
\draw (0,0) -- (2,0) -- (2,2) -- (0,2) -- (0,0); \draw (2/3,0) --
(2/3,2/3) -- (0,2/3); \draw (2,2/3) -- (4/3, 2/3) -- (4/3,2);; \draw
(2,4/3) -- (2/3,4/3) -- (2/3,2); \node at (1/3,1/3){$\si_1$}; \node
at (1,5/3){$\si_2$}; \node at (5/3,1){$\si_3$};
\end{tikzpicture}
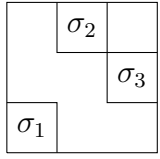
\caption{The diagram of $132[\si_1,\si_2,\si_3]$} \label{132} \efi

We wish to prove a third theorem of this type.  To do so, we will
need the notion inflation for permutations.  Given a permutation
$\pi=a_1\ldots a_n\in\fS_n$ and permutations $\si_1,\ldots,\si_n$,
the {\em inflation} of $\pi$ by the $\si_i$, written
$\pi[\si_1,\ldots,\si_n]$, is the permutation whose diagram is
obtained from the diagram of $\pi$ by replacing the dot $(i,a_i)$ by
a copy of $\si_i$ for $1\le i\le n$.  By way of example,
Figure~\ref{132} shows a schematic diagram of an inflation of the
form $132[\si_1,\si_2,\si_3]$.  More specifically,
$132[21,1,312]=216534$.

\bth Suppose $n=2^m-1$ for some  $m\ge0$.  Then
$$
[t^k]M_n(1,t)=A_{n,k}(1)=\case{1}{if $k=0$,}{\mbox{an even
integer}}{if $k\ge1$.}
$$
\eth \pf\ We have  $A_{n,0}(q)=1$ since $\si=12\ldots n$ is the only
permutation without descents and it avoids $321$.

If $k\ge1$, and $A_{n,k}(q)$ has an even number of terms then
$A_{n,k}(1)$ must be even because it is a symmetric polynomial by
Theorem~\ref{sym}. By the same token, if $A_{n,k}(q)$ has an odd
number of terms, then $A_{n,k}(1)$ has the same parity as its middle
term.  Now  consider $R_{180}$ acting on the elements of
$\cA_{nk/2}$ as in the previous proof.  Note that since $n$ is odd,
$k$ must be even.  Furthermore, this action partitions $\cA_{nk/2}$
into orbits of size one and two.  So it suffices to show that there
are an even number of fixed points.  If $\si$ is fixed then its
diagram must contain the center, $c$,  of the square since this is a
fixed point of $R_{180}$.  Also, the NW and SE quadrants of $\si$
with respect to $c$ must be empty, since otherwise they both must
contain dots (as one is taken to the other by $R_{180}$) and
together with $c$ this forms a $321$.   For the same reason, the SW
quadrant of $\si$ determines the NE one. Thus, the fixed points are
exactly the inflations of the form $\si=123[\tau,1,R_{180}(\tau)]$
where $\tau\in\Av_{2^{m-1}-1}(321)$ has $k/2$ descents.  By
induction on $m$ we have that  the number of such $\tau$, and hence the number of
such $\si$, is even. \hqedm

%
%

\section{A refinement of Theorem~\ref{I321:rr} using continued fractions}
\label{ptc}

We will now use a modification of the bijection in Section~\ref{ptm}
together with the theory of continued fractions to give a third
proof of Theorem~\ref{I321:rr}.  In fact, we will be able to keep
track of a third statistic on permutations $\si$, namely
$$
\fix\si=\text{the number of fixed points of $\si$.}
$$
So consider the following polynomial
$$
I_n(q,t,x)=\sum_{\si\in\Av_n(321)} q^{\inv\si} t^{\lrm\si}
x^{\fix\si}
$$
and the generating function
$$
{\mathfrak I}(q,t,x;z)=\sum_{n\ge0} I_n(q,t,x) z^n.
$$
It is worth noting that
 \beq\label{eq:exc-pol}
 I_n(q,t,x/t)=\sum_{\si\in\Av_n(321)} q^{\inv\si}
 t^{\exc\si}x^{\fix\si}
 \eeq
 where $\exc\si$ is the number of  excedances of $\si$ (i.e., the number of indices $i$ such that $\si(i)>i$).
 This follows from the following fact.

 \ble  \label{lem: exc vs lrm}
  Suppose $\si=a_1 a_2\cdots a_n\in\Av_n(321)$. Then $a_i$ is a left-right maximum if and only if $a_i\geq i$.
Consequently, $\lrm\si=\exc\si+\fix\si$.
 \ele

\prf If $a_i$ is a left-right maximum, then $a_i$ is greater
than the $i-1$ elements to its left in $\si$, so $a_i>i-1$.
Conversely, if $a_i$ is not a left-right maximum then  it is smaller than some element to its left in $\si$.  Also, by Lemma~\ref{321:structure}(a), it is smaller than all $n-i$ elements to its right in $\si$.
This implies that
$n-a_i\ge n-i+1$, whence $a_i\le i-1$.
 \hqedm

Continued fractions are very useful for enumerating weighted Motzkin
paths $M$.  If $s$ is a step of $M$ then we defined its {\em height}
to be
$$
\HT(s)=\text{the $y$-coordinate of the initial lattice point of $s$.}
$$
If $\HT(s)=i$ then we assign $s$ a weight $\wt s=u_i$, $d_i$, or
$l_i$ corresponding to $s$ being an  up, down, or level step,
respectively. Weight paths $M=s_1\ldots s_n$ and the set
$\cM_n$ of all such Motzkin paths by
$$
\wt M=\prod_{j=1}^n \wt s_j
$$
and
$$
\wt\cM_n=\sum_{M\in\cM_n} \wt M.
$$
For example, taking the path $M$ in Figure~\ref{fig:motzkin} would
give (ignoring the subscripts on the $L$'s) $\wt M=
u_0l_1^2u_1l_2d_2l_1d_1$.

In the sequel, we will use the following notation for continued
fractions \beq \label{cf} F=
\frac{a_1|}{|b_1}\pm\frac{a_2|}{|b_2}\pm\frac{a_3|}{|b_3}\pm\cdots
\quad = \quad \cfrac{a_1}{b_1\pm \cfrac{a_2}{b_2\pm
\cfrac{a_3}{b_3\pm\dotsb }}}. \eeq We can now state Flajolet's
classic result connecting continued fractions and weighted Motzkin
paths.

\bth[Flajolet~\cite{fla:cac}] \label{wt:thm} If $z$ is an
indeterminate then
 \beq \label{wt:eqn} \sum_{n\ge0}\wt \cM_n\  z^n
 =\frac{1|}{|1-l_0z}-\frac{u_0d_1z^2|}{|1-l_1z}-\frac{u_1d_2z^2|}{|1-l_2z}-\frac{u_2d_3z^2|}{|1-l_3z}-\cdots
\eeq is the generating function for weighted Motzkin paths.\hqed
\eth

We will also use the same classification of continued fractions
employed by Flajolet~\cite{fla:cac} which defines a Jacobi-type
continued fraction as a continued fraction of the form
\begin{align*}
J(z)&=\frac{1|}{|1-b_0 z}-\frac{ \la_1 z^2|}{|1-b_1 z}-\frac{\la_2
z^2|}{|1-b_2 z}
 -\frac{\la_3 z^2|}{|1-b_3z}-\cdots
\end{align*}
and a Stieljes-type continued fraction as a continued fraction of
the form
\begin{align*}
S(z)&=\frac{1|}{|1}-\frac{ \la_1 z|}{|1}-\frac{\la_2 z|}{|1}
 -\frac{\la_3 z|}{|1}-\cdots.
\end{align*}
Following the book of Jones and Thron~\cite{jt:cf}, a continued fraction of the form
\begin{align*}
T(z)&=\frac{1|}{|1-b_0 z}-\frac{ \la_1 z|}{|1-b_1 z}-\frac{\la_2
z|}{|1-b_2 z}
 -\frac{\la_3 z|}{|1-b_3z}-\cdots.
\end{align*}
will be called a Thron-type continued fraction.

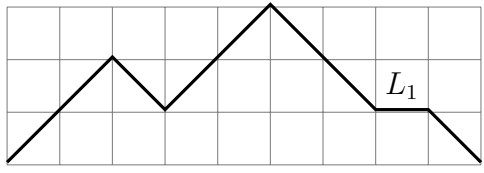
\begin{figure}
\begin{center}
\begin{tikzpicture}[scale=.7]
\draw[step=1,color=gray] (0,0) grid (9,3); \draw [very thick]
(0,0.05)--(1,1.05) --
(2,2.05)--(3,1.05)--(4,2.05)--(5,3.05)--(6,2.05)--(7,1.05)--(8,1.05)--(9,0.05);
\draw (7.5,1.5) node{$L_1$};
\end{tikzpicture}
\end{center}
\caption{The Motzkin path $\nu(\si)$ associated with
$\si=361782495$} \label{fig:motzkin2}
\end{figure}

In order to derive a continued fraction expansion for
$\fI(q,t,x;z)$, we will set up a bijection between $\Av_n(321)$ and
a subset of $\cM_n^{(2)}$. Call $M\in\cM_n^{(2)}$ {\em restricted} if
it has no $L_0$ steps at height~$0$.  Let $\cR_n$ be the set of such
paths.  In the following proof, we will use the same definitions and
notation as in Section~\ref{ptm}.
 \bth \label{fI:thm} The series $\fI(q,t,x;z)$ has Jacobi-type continued fraction expansion
$$
\fI(q,t,x;z)=\frac{1|}{|1-t x z}-\frac{ t qz^2|}{|1-(1+t)qz}-\frac{t
q^3 z^2|}{|1-(1+t)q^2z}
 -\frac{t q^5 z^2|}{|1-(1+t)q^3z}-\cdots.
$$
\eth \prf We define a bijection $\nu:\Av_n(321)\ra\cR_n$ in a way
similar to the bijection $\mu$ in the first proof of
Theorem~\ref{I321:rr}, but without the shift.
Specifically, given
$\si=a_1\ldots a_n\in\Av_n(321)$ with $\val\si=(v_1,\ldots,v_n)$ and
$\pos\si=(p_1,\ldots,p_n$), we let $\nu(\si)=M=s_1\ldots s_n$ where
$$
s_i=
\begin{cases}
U&\text{if $v_i=0$ and $p_i=1$,}\\
D&\text{if $v_i=1$ and $p_i=0$},\\
L_0&\text{if $v_i=p_i=0$,}\\
L_1&\text{if $v_i=p_i=1$.}
\end{cases}
$$
Continuing the example from the beginning of the paper,
$\si=361782495$ would be mapped to the path in
Figure~\ref{fig:motzkin2}.

We must show that $\nu$ is well defined in that $M\in\cR_n$.
Defining the inverse map and proving it is well defined is similar
and so left to the reader. The fact that $M$ is a Motzkin path
follows because, by Lemma~\ref{321:structure}(b), in every prefix of
$\pos\si$ the number of ones is at least as great as the number in
the corresponding prefix of $\val\si$, with equality for all of
$\si$. This forces similar inequalities and equality between the
number of up-steps and the number of down-steps in $M$.  Thus $M$
stays weakly above the $x$-axis and ends on it.

To see that $M$ has no $L_0$-steps on the $x$-axis, note first that
all steps before the first $U$-step (if any) must be of the form
$L_1$ because, if not, then the index $i$ of the first such $L_0$-step would contradict Lemma~\ref{321:structure}(b).  Also, any time $M$ returns to
the $x$-axis, it must be with $s_j=D$ for some $j$.  So the
corresponding prefixes of $\val\si$ and $\pos\si$ have the same
number of ones and this implies that $a_1\ldots a_j$ are
$1,\ldots,j$ in some order.  Now using an argument similar to the one just given, one sees that there can be no $L_0$-step before the next $U$-step.  This completes the proof that $\mu$ is well
defined.

We now claim that
\begin{eqnarray}
\label{fixsi}
\fix\si&=&\mbox{the number of $L_1$-steps at height $0$,}\\
\label{lrmsi2}
\lrm\si&=&\#U(M)+\#L_1(M),\\
\label{invsi2} \inv\si&=&\area M.
\end{eqnarray}

Let us prove the first equation.  If $s_j=L_1$ with $\HT(s_j)=0$
then, as in the proof that $\nu$ is well defined, $a_1\ldots
a_{j-1}$ are the numbers $1,\ldots,j-1$ in some order. Thus if
$v_j=p_j=1$ then both the position and value of $a_j$ correspond to
a left-right maximum. This forces $a_j=j$ and so we have a fixed
point.  Similar considerations show that every fixed point in a
$321$-avoiding permutation is a left-right maximum corresponding to
an $L_1$ step at height $0$.

Equation~\ree{lrmsi2} follows immediately from the fact that the
number of left-right maxima in $\si$ equals the number of ones in
$\pos\si$, and the corresponding steps in $M$ are of the form $U$ or
$L_1$.

For the final equality, first recall that all inversions of $\si$
are between a left-right maximum $m$ and a non-left-right maximum to
its right by  Lemma~\ref{321:structure}(a). So if $m$ is in position
$p$ then, because everything to its left is smaller, it creates
$m-p$ inversions.  Also, the maximum values and their positions in
$\si$ are given by $iv_i$ and $jp_j$, respectively, whenever
$v_i,p_j=1$.  Since $iv_i=jp_j=0$ whenever $v_i,p_j=0$ we have
$$
\inv\si=\sum_{v_i=1}iv_i-\sum_{p_j=1} jp_j=\sum_{i=1}^n (v_i-p_i)i.
$$

As far as the area, we start by noting that  $\area M=\sum_j \HT
s_j$.  Furthermore, $\HT(s_j)$ is just the difference between the
number of up-steps and down-steps preceding $s_j$. For any step
$s_i$, we have $p_i-v_i=1,-1,$ or $0$ corresponding to $s_i$ being
an up-, down-, or level-step, respectively. So $\HT(s_j)=\sum_{i<j}
(p_i-v_i)$. Combining expressions,  interchanging summations, and
using the fact that $\val\si$ and $\pos\si$ have the same number
of ones,  gives
$$
\area M=\sum_{j=1}^n\sum_{i<j} (p_i-v_i)=\sum_{i=1}^n
[(n-i)p_i-(n-i)v_i]=\sum_{i=1}^n (v_i-p_i)i.
$$
Comparing this expression with the one derived for $\inv\si$ in the
previous paragraph completes the proof of~\ree{invsi2}.

To finish the demonstration of the theorem, we just need to set the
weights in Theorem~\ref{wt:thm}  in light
of~\ree{fixsi}--\ree{invsi2}. The only level steps at height $0$ are
$L_1$ which contribute to both $\fix\si$ and $\lrm\si$.  So we let
$l_0=tx$.  At heights $h\ge1$ we have both $L_0$ and $L_1$ steps.
The former only contribute to $\inv$ by adding $h$, while the later
also increase the $\lrm$ statistic by one, so we have
$l_h=(1+t)q^h$.  Similar reasoning gives $u_h=tq^h$ and $d_h=q^h$
which, after plugging into equation~\ree{wt:eqn}, completes the
proof. \hqedm

 The following refinement of Theorem~\ref{I321:rr} is a simple consequence of the preceding result.

 \bth
For $n\ge1$,
$$
I_n(q,t,x)=txI_{n-1}(q,t,x)+\sum_{k=0}^{n-2} q^{k+1}
I_k(q,t,1)\left[ I_{n-1-k}(q,t,x)-t(x-1)I_{n-2-k}(q,t,x)\right].
$$
\eth

\prf We can derive from the continued fraction expansion of
$\fI(q,t,x;z)$ in Theorem~\ref{fI:thm} that
$$
\fI(q,t,x;z)=\cfrac{1}{1-t x z-\cfrac{t q z^2}{-qz
+\cfrac{1}{\fI(q,t,1;qz)}}}.
$$
After simplification, this leads to the functional equation
 \beq \label{eq:funcEqu}
  \fI(q,t,x;z)=1+txz\fI(q,t,x;z)+qz\fI(q,t,1;qz)\left[\fI(q,t,x;z)-1-tz(x-1)\fI(q,t,x;z)\right].
 \eeq
Extracting the coefficient of $z^n$ on both sides gives the desired
recursion. \hqedm

We note that the $q=1$ case of the functional
equation~\eqref{eq:funcEqu} is, by~\eqref{eq:exc-pol}, equivalent to
 equation~(1) in \cite{eli:fpe}. It can be explicitly solved as done in equation~(2) of the work just cited.

 One can simplify the continued fraction in Theorem~\ref{fI:thm} in the
case $x=1$.  The {\em $n$th convergent} of the continued
fraction~\ree{cf} is
$$
c_n(F)=\frac{a_1|}{|b_1}\pm\frac{a_2|}{|b_2}\pm\cdots\pm\frac{a_n|}{|b_n}.
$$
We wish to construct continued fractions, called the {\em even} and
{\em odd parts} of $F$ and denoted $F_e$ and $F_o$, such that
$c_{2n}(F)=c_n(F_e)$ and $c_{2n+1}(F)=c_n(F_o)$, respectively.
 The following theorem shows how to do this when $b_n=1$ for all $n$.
\bth[\cite{jt:cf}] \label{jt:thm} We have \bea
\frac{a_1|}{|1}+\frac{a_2|}{|1}+\frac{a_3|}{|1}+\cdots &=&
\frac{a_1|}{|1+a_2}-\frac{a_2a_3|}{|1+a_3+a_4}-\frac{a_4a_5|}{|1+a_5+a_6}-\frac{a_6a_7|}{|1+a_7+a_8}+\cdots\\[10pt]
&=&
a_1-\frac{a_1a_2|}{|1+a_2+a_3}-\frac{a_3a_4|}{|1+a_4+a_5}-\frac{a_5a_6|}{|1+a_6+a_7}
-\frac{a_7a_8|}{|1+a_8+a_9}-\cdots \eea as even and odd parts,
respectively,  of the first continued fraction.\hqed \eth

Using the concept of even and odd parts, we see that there is a
simple and well-known relationship between Stieljes- and Jacobi-type
continued fractions~\cite[p.\ 129]{jt:cf}.  More precisely, taking
the even part of a Stieljes-type continued fraction gives
$$
\frac{1|}{|1}-\frac{\la_1z|}{|1}-\frac{\la_2z|}{|1}-\cdots =
\frac{1|}{|1-\la_1z}-\frac{\la_1\la_2z^2|}{|1-(\la_2+\la_3)z}-\frac{\la_3\la_4z^2|}{|1-(\la_4+\la_5)z}-\frac{\la_5\la_6z^2|}{|1-(\la_6+\la_7)z}-\cdots
$$
where the right-hand side is of Jacobi type.  Combining
Theorem~\ref{fI:thm} when $x=1$ with the above equation, we get the
following result.

 \bco Set $\fI(q,t;z):=\fI(q,t,1;z)$. The generating function $\fI(q,t;z)$
has the Stieltjes-type continued fraction expansion \beq \label{sti}
\fI(q,t;z) =\frac{1|}{|1}-\frac{tz|}{|1}-\frac{qz|}{|1}
 -\frac{t qz|}{|1}-\frac{q^2 z|}{|1}-\frac{tq^2 z|}{|1}-\frac{q^3 z|}{|1}-\frac{tq^3 z|}{|1}
-\frac{q^4 z|}{|1}-\frac{tq^4 z|}{|1}-\cdots \eeq \eco

It is interesting to note that there is a second recursion for
$I_n(q,t)$ which follows from a result of Krattenthaler.
\bth[Krattenthaler~\cite{kra:prp}] We have
$$
\fI(q,t;z)=\frac{1|}{|1-(t-1)z}-\frac{z|}{|1-(tq-1)z}-\frac{z|}{|1-(tq^2-1)z}
 -\frac{z|}{|1-(tq^3-1)z}-\cdots
$$
as the Thron-type continued fraction expansion of $\fI(q,t;z)$.\hqed
\eth \bco\label{I321:rr2} For $n\ge1$,
$$
 I_n(q,t)=tI_{n-1}(q,t)+\sum_{k=0}^{n-2}I_k(q,t)I_{n-1-k}(q,qt).
$$
\eco \prf Simple manipulation of the continued fraction in
Krattenthaler's Theorem gives the functional equation
$$
 \fI(q,t;z)=1+(t-1)z \fI(q,t;z)+z \fI(q,t;z) \fI(q,qt;z).
$$
Taking the coefficient of $z^n$ on both sides of this equation
finishes the proof. \hqedm

It is worth noting that, {\it a priori}, it is not at all clear that the
recursions in Theorem~\ref{I321:rr} and the above corollary generate
the same sequence of polynomials.
The relationship between these two recursions can be interpreted in terms of the statistics sumpeaks and sumtunnels, introduced in Section~\ref{sss}, as follows.
Using the standard decomposition of non-empty Dyck paths as $P=NQER$, where $Q$ and $R$ are Dyck paths,
the generating function for $\cD_n$ with respect to the statistics $(\spea,\npea)$ satisfies
the recursion in Corollary~\ref{I321:rr2}. On the other hand, using the same decomposition, the generating function for $\cD_n$ with respect to the statistics $(\stun,n-\des)$
satisfies the recursion in Theorem~\ref{I321:rr}.
Thus, Theorem~\ref{thm:speastun} implies that the two recursions are equivalent.

%
%

\section{Refined sign-enumeration of $321$-avoiding permutations}
\label{rse}

Simion and Schmidt~\cite{ss:rp} considered the signed
enumeration of various permutation classes of the form
$\sum_{\si\in\Av_n(\pi)}(-1)^{\inv\si}$.  In this section we will
rederive their theorem for $\Av_n(321)$ using the results of the
previous section.  In addition, we will  provide a more refined signed
enumeration which also keeps track of the $\lrm$ statistic.  We
should note that Reifegerste~\cite{rei:rsb} also has a refinement
which takes into account the length of the longest increasing
subsequence of $\si$.

Let $C_n$ be the $n$th Catalan number and consider the generating
function $C(z)=\sum_{n\ge0} C_n z^n$. It is well known that $C(z)$
satisfies the functional equation $C(z)=1+zC(z)^2$. Rewriting this
as $C(z)=1/(1-z C(z))$ and iteratively substituting for $C(z)$, we
obtain the also well-known continued fraction
$$
C(z)=\frac{1|}{|1}-\frac{z|}{|1}-\frac{z|}{|1}-\frac{z|}{|1}-\frac{z|}{|1}-\frac{z|}{|1}-\cdots.
$$

Now plug $q=-1$ and $t=1$ into the continued fraction~\ree{sti} to
obtain
$$
\fI(-1,1;z)=\frac{1|}{|1}-\frac{z|}{|1}+\frac{z|}{|1}
 +\frac{z|}{|1}-\frac{z|}{|1}-\frac{z|}{|1}+\frac{z|}{|1}
 +\frac{z|}{|1}-\frac{z|}{|1}-\frac{z|}{|1}+\frac{z|}{|1}
 +\frac{z|}{|1}\cdots.
$$
Using Theorem~\ref{jt:thm} to extract the odd part of this expansion
gives
$$
\fI(-1,1;z)=1+\frac{z|}{|1}-\frac{z^2|}{|1}-\frac{z^2|}{|1}-\frac{z^2|}{|1}-\frac{z^2|}{|1}-\frac{z^2|}{|1}-\cdots.
$$
Comparing this to the continued fraction for $C(z)$, we see that
$$
\fI(-1,1;z)=1+zC(z^2).
$$
Taking the coefficient of $z^n$ on both sides yields the following
result. \bth[Simion and Schmidt~\cite{ss:rp}] For all $n\ge1$, we
have
$$
 I_{2n}(-1,1)=\sum_{\sigma\in\Av_{2n}(321)}(-1)^{\inv\si}=0\quad\text{and}\quad
 I_{2n+1}(-1,1)=\sum_{\sigma\in\Av_{2n+1}(321)}(-1)^{\inv\si}=C_{n}.
$$  \hqed \eth

Since our refined sign-enumeration will involve the parameter
$\lrm$, we recall (but will not use) the folklore result that the enumerating
polynomial of $\Av_n(321)$ according to the $\lrm$ statistic is the
$n$th {\em Narayana polynomial}, i.e.,
 $$
 I_{n}(1,t)=\sum_{\sigma\in\Av_{n}(321)}t^{\lrm\si}=\sum_{k=1}^n N_{n,k} t^k,
 $$
where the {\em Narayana number} $N_{n,k}$ is given by
$N_{n,k}=\frac{1}{n}{n\choose k}{n\choose k-1}$ for $n\ge k\ge1$.

\bth \label{thm:sign-enumeration} For all $n\geq 1$,
 \beq \label{eq:sign-enum}
  I_{n}(-1,t)=\sum_{\sigma\in\Av_{n}(321)}(-1)^{\inv\si}t^{\lrm\si}=\sum_{k=1}^{n}(-1)^{n-k} s_{n,k}
  t^{k}
 \eeq
where $s_{n,k}$ is defined for $n\ge k\ge1$ by
$$s_{n,k}={\flf{n-1}{2} \choose \flf{k-1}{2}} {\cef{n-1}{2} \choose\cef{k-1}{2}}.$$
 Moreover,
 \begin{eqnarray}
 I_{2n}(-1,t)&=&(t-1)I_{2n-1}(-1,t),\label{eq:rec1-sign}\\
 (n+1)I_{2n+1}(-1,t)&=&2((1+t^2)n-t)I_{2n-1}(-1,t)-(1-t^2)^2(n-1)I_{2n-3}(-1,t).\label{eq:rec2-sign}
 \end{eqnarray}
\eth

\prf Let $\fI(t;z)$ and $\fI_{odd}(t;z)$ be the power series defined
as
$$
\fI(t;z)=\sum_{n\ge0} I_{n}(-1,t)z^n\quad\text{and}\quad
\fI_{odd}(t;z)=\sum_{n\ge0} I_{2n+1}(-1,t)z^{n}.
$$
By equation~\ree{sti}, we have
\bea \fI(t;z) &=&
\frac{1|}{|1}-\frac{tz|}{|1}+\frac{z|}{|1}
 +\frac{tz|}{|1}-\frac{z|}{|1}-\frac{tz|}{|1}+\frac{z|}{|1}+\frac{tz|}{|1}
-\frac{z|}{|1}-\frac{tz|}{|1}+\frac{z|}{|1}+\frac{tz|}{|1}-\cdots
\\[10pt]
 &=&
\cfrac{1}{1- \cfrac{tz}{1+ \cfrac{z}{1+ \cfrac{tz}{1-z\fI(t;z)}}}}
\eea which, after simplification, leads to the functional equation
 $$
 (1+z-tz)z\,\fI(t;z)^2 -(1+2z+z^2-t^2z^2)\,\fI(t;z)+(1+z+tz)=0.
$$ Solving this quadratic equation, we obtain
 \beq\label{eq:GF-sign}
 \fI(t;z)=\frac{1+2z+(1-t^2)z^2-\sqrt{
 1-2(1+t^2)z^2+(1-t^2)^2z^4}}{2z(1+z-tz)}.
 \eeq
Noticing that
$\fI_{odd}(t;z)=(\fI(t;\sqrt{z})-\fI(t;-\sqrt{z}))/2\sqrt{z}$ and
using~\eqref{eq:GF-sign}, we obtain after a routine computation
 \beq\label{eq:GF-sign-odd}
 \fI_{odd}(t;z)=\frac{1-(1-t)^2z-\sqrt{1-2(1+t^2)z+(1-t^2)^2z^2}}{2z(1-(1-t)^2z)}.
 \eeq

 It follows from~\eqref{eq:GF-sign} that
$$
 (1+z-tz)\fI(t;z)+(1-z+tz)\fI(t;-z)=2.
$$
Extracting the coefficient of $z^{2n}$ on both sides of the last
equality, we obtain~\eqref{eq:rec1-sign}.
Using~\eqref{eq:GF-sign-odd}, it is easily checked  that $\fI_{odd}$
satisfies the differential equation
$$
 z\big(1-2(1+t^2)z+(1-t^2)^2z^2\big)\,\fI_{odd}'(t;z)+\big(1-2(1-t+t^2)z+(1-t^2)^2z^2)\big)\,\fI_{odd}(t;z)-t=0,
$$
 where $\fI_{odd}'(t;z)$ is the derivative with respect to $z$.
 Extracting the coefficient of $z^{n}$ on both sides of the last
equality, we obtain~\eqref{eq:rec2-sign}.\\

We now turn our attention to~\eqref{eq:sign-enum}. Clearly, we have
 \beq\label{eq:coeff-signenum-a}
 [t^{2k+1}]I_{2n+1}(-1,t)=[t^k z^n]\frac{\fI_{odd}(\sqrt{t};z)-\fI_{odd}(-\sqrt{t};z)}{2\sqrt{t}}
                         =[t^{k}][z^n] \frac{1}{\sqrt{1-2(1+t)z+(1-t)^2z^2}}
 \eeq
where the last equality follows from~\eqref{eq:GF-sign-odd}. Using
the Lagrange inversion formula, we can show that
 \beq
 [z^n]\frac{1}{\sqrt{1-2(1+t)z+(1-t)^2z^2}}=[x^n](1+(1+t)x+tx^2)^n.
 \eeq
Combining~\eqref{eq:coeff-signenum-a} with the above relation, we
obtain
 \beq\label{eq:coeff-signenum-for-a}
 [t^{2k+1}]I_{2n+1}(-1,t)=[x^n][t^{k}](1+(1+t)x+tx^2)^n=[x^n]\binom{n}{k}x^k(1+x)^n=\binom{n}{k}^2.
 \eeq
Similarly, we have
 \bea
 [t^{2k}]I_{2n+1}(-1,t)
  =[t^{k}z^n]\frac{\fI_{odd}(\sqrt{t};z)+\fI_{odd}(-\sqrt{t};z)}{2}
  =[t^{k}] [z^n]\frac{1}{2z}-\frac{1-z-tz}{2z\sqrt{1-2(1+t)z+(1-t)^2z^2}}.
 \eea
This, combined first with~\eqref{eq:coeff-signenum-a} and
then~\eqref{eq:coeff-signenum-for-a}, yields
$$
 [t^{2k}]I_{2n+1}(-1,t)=\frac{1}{2}\big([t^{2k+1}]I_{2n+1}(t,-1)+[t^{2k-1}]I_{2n+1}(t,-1)-[t^{2k+1}]I_{2n+3}(t,-1)\big)
  =-\binom{n}{k-1}\binom{n}{k}.
$$
This proves that~\eqref{eq:sign-enum} is true when $n$ is odd.
Combining this with~\ree{eq:rec1-sign} shows that the formula also holds when $n$ is even.
 \hqedm

To see why the previous result implies the one of Simion and Schmidt, plug $t=1$ into the equations for $I_{2n}(-1,t)$ and $I_{2n+1}(-1,t)$.  In the former case we immediately get $I_{2n}(-1,1)=0$ because of the factor of $t-1$ on the right.  In the latter, we get the equation $(n+1) I_{2n+1}(-1,1)=2(2n-1) I_{2n-1}(-1,1)$.  The fact that $I_{2n+1}(-1,1)=C_n$ now follows easily by induction.

Finally, it is interesting to note that the
numbers $s_{n,k}$ which arise in the signed enumeration of
$\Av_n(321)$ have a nice combinatorial interpretation.  Recall that
{\em symmetric} Dyck paths are those $P=s_1\ldots s_{2n}$ which are
the same read forwards as read backwards. The following result
appears in Sloane's Encyclopedia~\cite{slo:oei}: For $n\geq
k\geq1$, the number $s_{n,k}$ is equal to the number of symmetric
Dyck paths of semilength $n$ with $k$ peaks.

\end{document}